\documentclass{article}
\usepackage{amssymb,amsmath,latexsym}

\def\N{\mathbb{N}}
\def\Q{\mathbb{Q}}
\def\R{\mathbb{R}}
\def\Z{\mathbb{Z}}
\def\C{\mathbb{C}}
\def\B{{\overline B}}
\def\K{{\mathcal K}}
\def\c{\frak{c}}
\def\s{\frak{s}}
\def\m{{\bf m}}
\def\q{{\bf q}}
\def\r{{\bf r}}
\def\v{{\bf v}}
\def\cs{{\bf s}}
\newcommand\Res{\operatorname{Res}} 

\begin{document}
\setlength{\parindent}{0pt}
\setlength{\parskip}{0.4cm}

\newtheorem{theorem}{Theorem}
\newtheorem{lemma}[theorem]{Lemma}
\newtheorem{corollary}[theorem]{Corollary}
\newtheorem{definition}{Definition}

\begin{center}

\Large{\bf Dedekind cotangent sums}
\footnote{Appeared in \emph{Acta Arithmetica} {\bf 109}, no.~2 (2003), 109--130. \\ 
\emph{Keywords}: Dedekind sum, cotangent sum, Petersson-Knopp identity, lattice point enumeration in polyhedra. \\ 
\emph{2000 Mathematics Subject Classification.} {\em Primary} 11L03; {\em Secondary} 52C07.} 
\normalsize

{\sc Matthias Beck}

\end{center}

\footnotesize
{\it Abstract:} Let $ a, a_{1}, \dots , a_{d} $ be positive integers, $ m_{1} , \dots , m_{d} $ nonnegative
integers, and $ z_{1}, \dots , z_{d} $ complex numbers. We study expressions of the form
  \[ \sum_{ k \text{ mod } a } \ \prod_{ j=1 }^{ d } \cot^{ ( m_{j} ) } \pi \left( \frac{ k a_{j} }{ a } + z_{j} \right) \ . \]
Here the sum is taken over all $k$ mod $a$ for which the summand is not singular.
These {\bf Dedekind cotangent sums} generalize and unify various arithmetic sums introduced by
Dedekind, Rademacher, Apostol, Carlitz, Zagier, Berndt, Meyer, Sczech, and Dieter.
Generalized Dedekind sums appear in various areas such as analytic and algebraic number theory,
topology, algebraic and combinatorial geometry, and algorithmic complexity.
We prove reciprocity laws, Petersson-Knopp identities, and computability statements for the Dedekind cotangent sums.
\normalsize


\section{Introduction}
While studying the transformation properties of   
\[ \eta (z) := e^{ \pi i z / 12 } \prod_{ n \geq 1 } \left( 1 - e^{ 2 \pi i n z } \right) \ , \]
under $ \mbox{SL}_{2} ( \Z ) $, Dedekind, in the 1880's \cite{dedekind}, naturally arrived at the following arithmetic function:
Let $ ((x)) $ be the sawtooth function defined by
  \begin{equation}\label{doublebrackets} ((x)) := \left\{ \begin{array}{cl} \{ x \} - \frac{ 1 }{ 2 } & \mbox{ if } x \not\in \Z \\
                                                                            0                         & \mbox{ if } x \in \Z \ . \end{array} \right. \end{equation}
Here $ \{ x \} = x - [x] $ denotes the fractional part of $x$. 
For $a, b \in \N := \left\{ n \in \Z : n > 0 \right\} $, we define the {\bf Dedekind sum} as
  \begin{equation}\label{dedsaw} \s (a,b) := \sum_{ k \mbox{ \rm \scriptsize mod } b } \left( \left( \frac{ ka }{ b } \right) \right) \left( \left( \frac{ k }{ b } \right) \right) \ . \end{equation}

The Dedekind sums and their generalizations have since intrigued mathematicians from various areas
such as analytic \cite{dedekind,dieter,almkvist} and algebraic number theory \cite{meyer,solomon}, topology \cite{hirz,zagier,sczech}, 
algebraic \cite{pommersheim,brion} and combinatorial geometry \cite{mordell,diaz}, and algorithmic complexity \cite{knuth}. 

By means of the discrete Fourier series of the sawtooth function (see, for example, \cite[p.~14]{grosswald}),
  \begin{equation}\label{cotfourier} \left( \left( \frac{ n }{ p } \right) \right)  = \frac{ i }{ 2p } \sum_{ k = 1 }^{ p-1 } \cot \left( \frac{ \pi k }{ p } \right) e^{ 2 \pi i k n / p } \ , \end{equation}
it is not hard to write the Dedekind sum in terms of cotangents:
  \begin{equation}\label{dedcot} \s(a,b) = \frac{ 1 }{ 4b } \sum_{ k=1 }^{ b-1 } \cot \frac{ \pi k a }{ b } \cot \frac{ \pi k }{ b }\ .  \end{equation}
Starting with these two representations (\ref{dedsaw}) and (\ref{dedcot}) of $ \s (a, b) $, various
generalizations of the Dedekind sum have been introduced. 
This paper constitutes an attempt to unify generalizations of the Dedekind sum in its `cotangent 
representation' (\ref{dedcot}). Through the discrete Fourier series (\ref{cotfourier}), this includes most 
generalizations of the `sawtooth representation' (\ref{dedsaw}) of the Dedekind sum. 

Let $ \cot^{(m)} $ denote the $m$'th derivative of the cotangent function. Our generalized Dedekind sum 
is introduced as follows: 
\begin{definition}\label{defi} For $ a_{0}, \dots , a_{d} \in \N $, $ m_{0} , \dots , m_{d} \in \N_{0} := \N \cup \{ 0 \} $, $ z_{0}, \dots , z_{d} \in \C $,
we define the {\bf Dedekind cotangent sum} as
  \[ \c \left( \begin{array}{c|ccc} a_{0} & a_{1} & \cdots & a_{d} \\
                                    m_{0} & m_{1} & \cdots & m_{d} \\
                                    z_{0} & z_{1} & \cdots & z_{d} \end{array} \right)
        := \frac{1}{ a_{0}^{ m_{0} + 1 } } \sum_{ k \mbox{ \rm \scriptsize mod } a_{0} } \ \prod_{ j=1 }^{ d } \cot^{ ( m_{j} ) } \pi \left( a_{j} \frac{ k + z_{ 0 } }{ a_{ 0 } } - z_{j} \right) \ , \]
where the sum is taken over all $k$ {\rm mod} $ a_{0} $ for which the summand is not singular.
\end{definition}
The Dedekind cotangent sums include as special cases various generalized Dedekind sums 
introduced by Rademacher \cite{rademacher}, Apostol \cite{apostol}, Carlitz \cite{carlitzbernoulli}, Zagier \cite{zagier}, 
Berndt \cite{berndt}, Meyer, Sczech \cite{sczech}, and Dieter \cite{dietercot}. 
Section \ref{history} contains the definitions of these previously defined sums and the connections to our new definition. 

The most fundamental and important theorems for any of the generalized De\-de\-kind sums are the
{\it reciprocity laws}: an appropriate sum of generalized Dedekind sums (usually permuting the
arguments in a cyclic fashion) gives a simple rational expression. 
The famous reciprocity law for the classical Dedekind sum is as old as the sum itself: 
\begin{theorem}[Dedekind]\label{dedreclaw} If $a,b \in \N$ are relatively prime then
  \[ \s (a,b) + \s (b,a) = - \frac{1}{4} + \frac{1}{12} \left( \frac{ a }{ b } + \frac{ 1 }{ ab } + \frac{ b }{ a } \right) \ . \]
\end{theorem}
The following reciprocity law for the Dedekind cotangent sums will be proved in section \ref{reclaw}: 
\begin{theorem}\label{main} Let $ a_{0} , \dots , a_{d} \in \N $, $ m_{0} , \dots m_{d} \in \N_{0} $, $
z_{0} , \dots , z_{d} \in \C $. If for all distinct $ i, j \in \{ 0 , \dots , d \} $ and all $ m,n \in \Z $,
 \[ \frac{ m + z_{i} }{ a_{i} } - \frac{ n + z_{j} }{ a_{j} } \not\in \Z \ , \] 
then
  \begin{eqnarray*} &\mbox{}& \sum_{ n=0 }^{ d } (-1)^{ m_{n} } m_{n} ! \sum_{ { l_{0} , \dots , \widehat{ l_{n} } , \dots , l_{d} \geq 0 } \atop { l_{0} + \dots + \widehat{ l_{n} } + \dots + l_{d} = m_{n} } } \frac{ a_{0}^{ l_{0} } \cdots \widehat{ a_{n}^{ l_{n} } } \cdots a_{d}^{ l_{d} } }{ l_{0}! \cdots \widehat{ l_{n}! } \cdots l_{d}! } \ 
                                \c \left( \begin{array}{c|ccccc} a_{n} & a_{0} & \cdots & \widehat{ a_{n} } & \cdots & a_{d} \\
                                                                 m_{n} & m_{0} + l_{0} & \cdots & \widehat{ m_{n} + l_{n} } & \cdots & m_{d} + l_{d} \\
                                                                 z_{n} & z_{0} & \cdots & \widehat{ z_{n} } & \cdots & z_{d} \end{array} \right) \\
                    &\mbox{}& = \left\{ \begin{array}{ll} (-1)^{ d/2 } & \mbox{ if all } m_{k} = 0 \mbox{ and } d \mbox{ is even, } \\
                                                          0            & \mbox{ otherwise. } \end{array} \right. \end{eqnarray*}
\end{theorem}
As usual, $ \widehat{ x_{n} } $ means we omit the term $ x_{n} $.

In section \ref{var}, we obtain as corollaries the reciprocity laws for Dedekind's, Apostol's, 
Meyer-Sczech's, Dieter's and Zagier's generalized Dedekind sums. The purpose here is not so much 
a rederivation of old theorems but rather to show a common thread to all of them. 

Another basic identity on the classical Dedekind sum is about a hundred years younger than Dedekind's reciprocity law \cite{knopp}: 
\begin{theorem}[Petersson-Knopp]\label{petknopp} Suppose $a,b \in \N$ are relatively prime. Then
  \begin{equation} \sum_{ d|n } \ \sum_{ k \mbox{ \rm \scriptsize mod } d } \s \left( \frac{ n }{ d } b + ka , ad \right) = \sigma (n) \ \s ( b,a ) \ . \end{equation}
Here $ \sigma (n) $ denotes the sum of the positive divisors of $n$.
\end{theorem}
This result has also been extended to certain generalized Dedekind sums \cite{vu,rosen,zheng}. 
The respective identity for the cotangent sum follows from a much more general theorem, stated 
and proved in section \ref{petkn}: 
\begin{theorem}\label{petknoppcotangentsum} For $ n, a_{0}, \dots , a_{d} \in \N , m_{0}, \dots , m_{d} \in \N_{0} $,
  \begin{eqnarray*} &\mbox{}& \sum_{ b | n } b^{ m_{0} + 1 - m_{1} - \dots - m_{d} - d } \sum_{ r_{1} , \dots , r_{d} \mbox{ \rm \scriptsize mod } b }
                                 \c \left( \begin{array}{c|ccc} a_{0} b & \frac{ n }{ b } a_{1} + r_{1} a_{0} & \cdots & \frac{ n }{ b } a_{d} + r_{d} a_{0} \\
                                                                m_{0} & m_{1} & \cdots & m_{d} \\
                                                                    0 & 0 & \cdots & 0 \end{array} \right) \\
                    &\mbox{}& \qquad = \ n \ \sigma_{ - m_{1} - \dots - m_{d} - 1 } (n) \
                                 \c \left( \begin{array}{c|ccc} a_{0} & a_{1} & \cdots & a_{d} \\
                                                                m_{0} & m_{1} & \cdots & m_{d} \\
                                                                    0 & 0 & \cdots & 0 \end{array} \right) \ . \end{eqnarray*}
Here $ \sigma_{m} (n) := \sum_{ d|n } d^{m} $.
\end{theorem}
Finally, we show in section \ref{comp} a computability result:  
\begin{theorem}\label{compcotangent} The Dedekind cotangent sum \ 
   $ \c \left( \begin{array}{c|ccc} a_{0} & a_{1} & \cdots & a_{d} \\
                                    m_{0} & m_{1} & \cdots & m_{d} \\
                                    z_{0} & z_{1} & \cdots & z_{d} \end{array} \right) $ 
is polynomial-time computable in the input size of $ a_{0} , \dots , a_{d} $.
\end{theorem}
Computability for any of the generalized Dedekind sums involving two integer arguments follows easily 
from the two-term reciprocity laws. However, the reciprocity laws for `higher-dimensional' analogues 
do not imply nice computability properties. 


\section{Various Dedekind sums}\label{history} 
In this section we will give an overview of previously defined generalizations of the Dedekind sum. 
We do not claim any completeness but hope to give a good picture of what has been introduced in the past. 

The sawtooth function $((x))$ defined in (\ref{doublebrackets}) is the first Bernoulli function $ \B_{1} (x) $, defined to be zero at the
integers. We used the slightly antiquated function $((x))$ partly for historical reasons, partly because it has a simpler discrete Fourier series. 
The {\bf Bernoulli polynomials} $ B_{k} (x) $ are defined through 
  \begin{equation}\label{berdef} \frac{z e^{xz } }{e^{z } - 1 } = \sum_{k \geq 0 } \frac{B_{k} (x) }{k! } z^{k} \ . \end{equation}
The first few of them are
  \begin{eqnarray*} &\mbox{}& B_{1} (x) = x - \frac{1}{2} \\
                    &\mbox{}& B_{2} (x) = x^{2} - x + \frac{1}{6} \\
                    &\mbox{}& B_{3} (x) = x^{3} - \frac{3}{2} x^{2} + \frac{1}{2} x \\
                    &\mbox{}& B_{4} (x) = x^{4} - 2 x^{3} + x^{2} - \frac{1}{30} \\
                    &\mbox{}& B_{5} (x) = x^{5} - \frac{5}{2} x^{4} + \frac{5}{3} x^{3} - \frac{1}{6} x \\
                    &\mbox{}& B_{6} (x) = x^{6} - 3 x^{5} + \frac{5}{2} x^{4} + \frac{1}{2} x^{2} + \frac{1}{42} \\
                    &\mbox{}& B_{7} (x) = x^{7} - \frac{7}{2} x^{6} + \frac{7}{2} x^{5} + \frac{7}{6} x^{3} + \frac{1}{6} x \ . \end{eqnarray*} 
%
The {\bf Bernoulli numbers} are $ B_{k} := B_{k} (0) $. The {\bf Bernoulli functions} $ \B_{k} (x) $ are the
periodized Bernoulli polynomials:
  \[\B_{k} (x) := B_{k} ( \{ x \} ) \ .  \]
Apostol \cite{apostol} replaced one of the sawtooth functions in (\ref{dedsaw}) by an arbitrary Bernoulli function:
  \begin{equation}\label{apo} \sum_{ k \mbox{ \rm \scriptsize mod } b } \frac{k}{b} \ \B_{n} \left( \frac{ka }{b } \right) \ .  \end{equation}
Apostol's idea was generalized by Carlitz \cite{carlitzbernoulli} and Mikol\'as \cite{mikolas} 
to what we would like to call the {\bf Dedekind Bernoulli sum}, defined for $ a, b, c, m, n \in \N $ 
as 
  \[ \s_{ m,n } ( a; b, c ) := \sum_{ k \mbox{ \rm \scriptsize mod } b } \B_{m} \left( \frac{ kb }{ a } \right) \B_{n} \left( \frac{ kc }{ a } \right) \ . \] 
Another way of generalizing (\ref{dedsaw}) is to shift the argument of the sawtooth functions.
This was introduced by Meyer \cite{meyer} and Dieter \cite{dieter}, and brought to a solid
ground by Rademacher \cite{rademacher}: 
For $ a,b \in \N $, $ x,y \in \R $, the {\bf Dedekind-Rademacher sum} is defined by
  \[ \s (a,b;x,y) := \sum_{ k \mbox{ \rm \scriptsize mod } b } \left( \left( a \frac{ k+y }{ b } - x \right) \right) \left( \left( \frac{ k+y }{ b } \right) \right) \ . \]
Note that there is no loss in restricting $x$ and $y$ to $ 0 \leq x, y < 1 $.

The ideas of Apostol and Rademacher can also be combined: Tak\'acs \cite{takacs} introduced a shift in Apostol's sum (\ref{apo}):
  \[ \sum_{ k \mbox{ \rm \scriptsize mod } b } \B_{1} \left( \frac{k+y}{b} \right) \B_{n} \left( a \frac{k+y}{b } - x \right) \ . \]
This was further generalized by Halbritter \cite{halbritter} and later by Hall, Wilson, and Zagier \cite{hall} 
to the {\bf generalized Dedekind-Ra\-de\-ma\-cher sum}, defined for $ a,b,c, m, n \in \N $, $ x,y,z \in \R $, 
by 
  \[ \s_{ m,n } \left( \begin{array}{c|l} a & b \ c \\
                                          x & y \ z \end{array} \right) :=
     \sum_{ k \mbox{ \rm \scriptsize mod } a } \B_{m} \left( b \frac{k+x}{a} - y \right) \B_{n} \left( c \frac{k+x}{a} - z \right) \ . \] 
On the other hand, we can start with the cotangent representation (\ref{dedcot}) of the Dedekind
sum to arrive at generalizations. The equivalent of the Dedekind-Rademacher sum in terms of
cotangents was first defined by Meyer and Sczech \cite{sczech}, motivated by the appearance of the
classical Dedekind sum in topology \cite{hirz}. The version we state here was introduced by Dieter \cite{dietercot}: 
For $ a,b,c \in \N $, $ x,y,z \in \R $, the {\bf cotangent sum} is defined by
  \[ \c (a,b,c;x,y,z) := \frac{1}{c} \sum_{ k \mbox{ \rm \scriptsize mod } c } \cot \pi \left( a \frac{ k+z }{ c } - x \right) \cot \pi \left( b \frac{ k+z }{ c } - y \right) \ . \]
Here the sum is taken over all $k$ {\rm mod} $b$ for which the summand is not singular.
Dieter remarked in \cite{dietercot} that the cotangent sums include as special cases various
modified Dedekind sums introduced by Berndt \cite{berndt,berndt2}. Most of them
are inspired by the transformation properties of the logarithm of the classical theta-function.
We list them here using Berndt's notation; throughout, $a$ and $b$ denote relatively prime positive
integers.
In the first sum, $ \alpha, \beta \in \N $, and $ a^{-1} $ is defined through $ a^{-1} a \equiv 1 $ mod $b$.
  \begin{eqnarray*} &\mbox{}& s_{ \alpha , \beta } (a,b) := \sum_{ k=1 }^{ ab-1 } \exp \left( 2 \pi i \left( \frac{ k \alpha }{ a } + \frac{ k \beta }{ b } \right) \right) \left( \left( \frac{ k }{ ab } \right) \right) \left( \left( \frac{ k a^{-1} }{ b } \right) \right) \\
                    &\mbox{}& S ( a, b ) := \sum_{ k=1 }^{ b-1 } (-1)^{ k + 1 + [ ak/b ] } \\
                    &\mbox{}& s_{1} ( a , b ) := \sum_{ k=1 }^{ b } (-1)^{ [ ak/b ] } \left( \left( \frac{ k }{ b } \right) \right)  \\
                    &\mbox{}& s_{2} ( a , b ) := \sum_{ k=1 }^{ b } (-1)^{ k } \left( \left( \frac{ k }{ b } \right) \right) \left( \left( \frac{ ka }{ b } \right) \right)  \\
                    &\mbox{}& s_{3} ( a , b ) := \sum_{ k=1 }^{ b } (-1)^{ k } \left( \left( \frac{ ka }{ b } \right) \right)  \\
                    &\mbox{}& s_{4} ( a , b ) := \sum_{ k=1 }^{ b-1 } (-1)^{ [ ak/b ] } \\
                    &\mbox{}& s_{5} ( a , b ) := \sum_{ k=1 }^{ b } (-1)^{ k + [ ak/b ] } \left( \left( \frac{ k }{ b } \right) \right) \ . \end{eqnarray*}
Yet another generalization of (\ref{dedcot}) we mention here is due to Zagier \cite{zagier}. 
Generalizing the topological considerations in \cite{hirz} to arbitrary dimensions, he arrived naturally at the following expression:
Let $ a_{1} , \dots , a_{d} \in \N $ be relatively prime to $ a_{0} \in \N $. Define the {\bf higher-dimensional Dedekind sum} as
  \[ \s ( a_{0} ; a_{1} , \dots , a_{d} ) := \frac{ (-1)^{d/2} }{ a_{0} } \sum_{ k=1 }^{ a_{0} - 1 } \cot \frac{ \pi k a_{1} }{ a_{0} } \cdots \cot \frac{ \pi k a_{d} }{ a_{0} } \ . \]
We note that this sum vanishes if $d$ is odd, since the cotangent is an odd function.
Berndt noticed in \cite{berndt} that a version of the higher-dimensional Dedekind sum was already
introduced by Carlitz \cite{carlitz} via sawtooth functions.

Our Dedekind cotangent sums (Definition \ref{defi}) combine the cotangent and higher-dimensional Dedekind 
sum. This could have been done by only introducing a shift of the variable in each cotangent of the 
higher-dimensional Dedekind sum. The reason for additionally introducing cotangent derivatives is twofold: 
first, they appear in lattice point enumeration formulas for polyhedra \cite{diaz}. 
We will make use of this fact in section \ref{comp} about the computability of the Dedekind cotangent sums.
Second, the cotangent derivatives are essentially the discrete Fourier transforms of the Bernoulli
functions. That is, our definition can be seen as the higher-dimensional `cotangent equivalent' to Apostol's
Dedekind Bernoulli sum.
In fact, the Dedekind cotangent sums include as special cases all generalized Dedekind sums mentioned in this 
section with exception of the Dedekind-Rademacher sum and its generalization by Hall, Wilson, and Zagier. 
Technically, these sums could be treated in the same manner. What makes a subtle difference is the shift 
by a real number in the argument of the Bernoulli functions: to handle these, we would have to work with 
the discrete Fouries series of the shifted Bernoulli functions, which turn out to be much less practical than 
the ones without a shift. 

To be more precise, we state the discrete Fourier series of the Bernoulli functions next. In analogy to 
(\ref{cotfourier}), we have the straightforward 
\begin{lemma}\label{fourier} For $ m \geq 2 $,
  \[ \B_{m} \left( \frac{ n }{ p } \right) = \frac{ B_{m} }{ (-p)^{ m } } + m \left( \frac{ i }{ 2p } \right)^{ m } \ \sum_{ k = 1 }^{ p-1 } \cot^{ (m-1) } \left( \frac{ \pi k }{ p } \right) e^{ 2 \pi k n / p } \ . \] 
\hfill {} $\Box$
\end{lemma}

These discrete Fourier expansions can be used, for example, to rewrite the Dedekind Bernoulli sums in terms of the Dedekind cotangent sums:
\begin{corollary}\label{apocot} If $ a, b, c \in \N $ are pairwise relatively prime and $ m, n \geq 2 $ are integers with the same parity then
  \begin{align*} \s_{ m,n } ( a; b, c ) &\stackrel{ \mbox{\rm \scriptsize def} }{ = } \sum_{ k \mbox{ \rm \scriptsize mod } a } \B_{m} \left( \frac{ kb }{ a } \right) \B_{n} \left( \frac{ kc }{ a } \right) \\
                                        &\, = mn \frac{ (-1)^{ (m-n)/2 } }{ 2^{ m+n } a^{ m+n-1 } } \sum_{ k=1 }^{ a-1 } \cot^{ (m-1) } \left( \frac{ \pi k c }{ a } \right) \cot^{ (n-1) } \left( \frac{ \pi k b }{ a } \right) + \frac{ B_m B_n }{ a^{ m+n-1 } } \\
                                        &\stackrel{ \mbox{\rm \scriptsize def} }{ = } mn \frac{ (-1)^{ (m-n)/2 } }{ 2^{ m+n } }
                                           \ \c \left( \begin{array}{c|cc} a & b & c \\
                                                                       m+n-2 & n-1 & m-1 \\
                                                                           0 & 0 & 0 \end{array} \right)
                                           + \frac{ B_m B_n }{ a^{ m+n-1 } } \ . \end{align*}
\end{corollary}
We note that the parity assumption on $m$ and $n$ is no restriction, since the sums vanish if $ m+n $ is
odd.
It is worth mentioning that a close relative of these sums, namely,
  \[ \sum_{ k=1 }^{ a-1 } \B_{m} \left( \frac{ k }{ a } \right) \left( \left( \frac{ kb }{ a } \right) \right) = m \frac{ (-1)^{ (m-1)/2 } }{ 2^{ m+1 } a^{m} } \sum_{ k=1 }^{ a-1 } \cot \left( \frac{ \pi k }{ a } \right) \cot^{(m-1)} \left( \frac{ \pi kb }{ a } \right) \]
appears naturally in the study of plane partition enumeration \cite{almkvist}. 

{\it Proof of Corollary} \ref{apocot}: By Lemma \ref{fourier},
  \begin{align} \s_{ m,n } ( a; b, c ) &= \sum_{ k \mbox{ \rm \scriptsize mod } a } \B_{m} \left( \frac{ kb }{ a } \right) \B_{n} \left( \frac{ kc }{ a } \right) \nonumber \\
                                       &= mn \left( \frac{i}{2a} \right)^{ m+n } \sum_{ k \mbox{ \rm \scriptsize mod } a } \ \sum_{ j,l = 1 }^{ a-1 } \cot^{ (m-1) } \left( \frac{ \pi j }{ a } \right) \cot^{ (n-1) } \left( \frac{ \pi l }{ a } \right) e^{ 2 \pi i k (jb+lc) / a } \nonumber \\
                                       &\quad + \ m \left( \frac{i}{2a} \right)^{ m } \frac{ B_{n} }{ (-a)^{n} } \sum_{ k \mbox{ \rm \scriptsize mod } a } \ \sum_{ j = 1 }^{ a-1 } \cot^{ (m-1) } \left( \frac{ \pi j }{ a } \right) e^{ 2 \pi i k jb / a } \label{mess} \\
                                       &\quad + \ n \left( \frac{i}{2a} \right)^{ n } \frac{ B_{m} }{ (-a)^{m} } \sum_{ k \mbox{ \rm \scriptsize mod } a } \ \sum_{ l = 1 }^{ a-1 } \cot^{ (n-1) } \left( \frac{ \pi l }{ a } \right) e^{ 2 \pi i k lc / a } \nonumber + \ a \frac{ B_{m} B_{n} }{ (-a)^{m+n} } \ . \nonumber \end{align}
Now we use the fact that $m+n$ is even, $ a,b,c $ are pairwise relatively prime, and
  \[ \sum_{ k \mbox{ \rm \scriptsize mod } a } e^{ 2 \pi i n k / a } = \left\{ \begin{array}{cl} a & \mbox{ if } a|n \\
                                                                                                 0 & \mbox{ else. } \end{array} \right. \]
Hence the first sum in (\ref{mess}) vanishes unless $ jb+lc $ is divisible by $a$, that is,
$ l \equiv - j b c^{-1} $ mod $a$, where $ c^{-1} c \equiv 1 $ mod $a$.
The second and third sum in (\ref{mess}) disappear completely, and we are left with
  \begin{align*} \s_{ m,n } ( a; b, c ) &= mn \frac{ (-1)^{ (m+n)/2 } }{ (2a)^{ m+n } } \ a \sum_{ j = 1 }^{ a-1 } \cot^{ (m-1) } \left( \frac{ \pi j }{ a } \right) \cot^{ (n-1) } \left( - \frac{ \pi j b c^{-1} }{ a } \right) + \frac{ B_{m} B_{n} }{ a^{m+n-1} } \\
                                        &= mn \frac{ (-1)^{ (m+n)/2 } }{ 2^{ m+n } a^{ m+n-1 } } \ (-1)^{n} \sum_{ j = 1 }^{ a-1 } \cot^{ (m-1) } \left( \frac{ \pi j c }{ a } \right) \cot^{ (n-1) } \left( \frac{ \pi j b }{ a } \right) + \frac{ B_{m} B_{n} }{ a^{m+n-1} } \ . \end{align*} 
 \hfill {} $\Box$


\section{Proof of the reciprocity law}\label{reclaw} 
{\it Proof of Theorem} \ref{main}. Consider the function
  \[ f(z) = \prod_{ j=0 }^{ d } \cot^{ ( m_{j} ) } \pi \left( a_{j} z - z_{j} \right) \ . \]
We integrate $f$ along the simple rectangular path
  \[ \gamma = [ x+iy, x-iy, x+1-iy, x+1+iy, x+iy ] \ , \]
where $x$ and $y$ are chosen such that $ \gamma $ does not pass through any pole of $f$, and all poles $ z_{p} $ of $f$ have imaginary
part $ | \mbox{Im} ( z_{p} ) | < y $. By the periodicity of the cotangent, the contributions of the two vertical segments of $ \gamma $
cancel each other. By definition of the cotangent,
  \[ \lim_{ y \to \infty } \cot ( x \pm iy ) = \mp i \ , \]
and therefore also
  \[ \lim_{ y \to \infty } \cot^{ (m) } ( x \pm iy ) = 0 \]
for $ m>0 $. Hence if any of the $ m_{j} > 0 $,
  \[ \int_{ \gamma } f(z) \ dz = 0 \ . \]
If all $ m_{j} = 0 $, we obtain
  \[ \int_{ \gamma } f(z) \ dz = i^{ d+1 } - (-i)^{ d+1 } = i^{ d+1 } \left( 1 + (-1)^{d} \right) \ . \]
This can be summarized in
  \[ \frac{1}{2i} \int_{ \gamma } f(z) \ dz = \left\{ \begin{array}{ll} i^{ d } & \mbox{ if all } m_{j} = 0 \mbox{ and } d \mbox{ is even } \\
                                                                        0       & \mbox{ otherwise, }
\end{array} \right. \]
or, by means of the residue theorem,
  \begin{equation}\label{residues} \pi \sum_{ z_{p} } \Res \left( f(z) , z_{p} \right) = \left\{ \begin{array}{ll} i^{ d } & \mbox{ if all } m_{j} = 0 \mbox{ and } d \mbox{ is even } \\
                                                                                                                           0       & \mbox{ otherwise. } \end{array} \right. \end{equation}
Here the sum ranges over all poles $ z_{p} $ inside $ \gamma $. It remains to compute their residues.
By assumption, $f$ has only simple poles. We will compute the residue at $ z_{p} = \frac{ k + z_{0} }{ a_{0} } $,
$ k \in \Z $, the other residues are completely equivalent. We use the Laurent expansion of the cotangent
  \begin{equation}\label{laurent} \cot \pi \left( a_{0} z - z_{0} \right) = \frac{ 1 }{ \pi a_{0} } \left( z - \frac{ k + z_{0} }{ a_{0} } \right)^{-1} + \mbox{ analytic part } \ , \end{equation}
and, more generally,
  \[ \cot^{ ( m_{0} ) } \pi \left( a_{0} z - z_{0} \right) = \frac{ (-1)^{ m_{0} } m_{0}! }{ ( \pi a_{0} )^{ m_{0} + 1 } } \left( z - \frac{ k + z_{0} }{ a_{0} } \right)^{-( m_{0} + 1)} + \mbox{ analytic part } \ . \]
The other cotangents are analytic at this pole: for $ j > 0 $,
  \begin{equation}\label{taylor} \cot^{ ( m_{j} ) } \pi \left( a_{j} z - z_{j} \right) = \sum_{ l_{j} \geq 0 } \frac{ ( \pi a_{j} )^{ l_{j} } }{ l_{j}! } \cot^{ ( m_{j} + l_{j} ) } \pi \left( a_{j} \frac{ k + z_{0} }{ a_{0} } - z_{j} \right) \left( z - \frac{ k + z_{0} }{ a_{0} } \right)^{ l_{j} } \ . \end{equation}
Hence
  \begin{eqnarray*} &\mbox{}& \Res \left( f(z) , z = \frac{ k + z_{0} }{ a_{0} } \right) = \frac{ (-1)^{ m_{0} } m_{0}! }{ \pi a_{0}^{ m_{0} + 1 } } \sum_{ { l_{1} , \dots , l_{d} \geq 0 } \atop { l_{1} + \dots + l_{d} = m_{0} } } \prod_{ j=1 }^{ d } \frac{ a_{j}^{ l_{j} } }{ l_{j}! } \cot^{ ( m_{j} + l_{j} ) } \pi \left( a_{j} \frac{ k + z_{0} }{ a_{0} } - z_{j} \right) \ . \end{eqnarray*}
Since $ \gamma $ has horizontal width 1, we have $ a_{0} $ poles of the form $ \frac{ k + z_{0} }{ a_{0} } $ inside $ \gamma $,
where $k$ runs through a complete set of residues modulo $ a_{0} $. This gives, by definition of the Dedekind cotangent sum,
  \[ \sum_{ k \mbox{ \rm \scriptsize mod } a_{0} } \Res \left( f(z) , z = \frac{ k + z_{0} }{ a_{0} } \right) = 
       \frac{1}{\pi} (-1)^{ m_{0} } m_{0}! \sum_{ { l_{1} , \dots , l_{d} \geq 0 } \atop { l_{1} + \dots + l_{d} = m_{0} } } \frac{ a_{1}^{ l_{1} } \cdots a_{d}^{ l_{d} } }{ l_{1}! \cdots l_{d}! }
       \ \c \left( \begin{array}{c|ccc} a_{0} & a_{1} & \cdots & a_{d} \\
                                        m_{0} & m_{1} + l_{1} & \cdots & m_{d} + l_{d} \\
                                        z_{0} & z_{1} & \cdots & z_{d} \end{array} \right) \ . \] 
The other residues are computed in the same way, and give with (\ref{residues}) the statement.
\hfill {} $\Box$


\section{Variations on a theme}\label{var}
The conditions on the parameters appearing in Theorem \ref{main} are not crucial; however,
without them the theorem would not be as easy to state. Indeed, the 
conditions simply ensure that all poles of the function $f$ used in the proof are simple.
We will now state some special cases of Theorem \ref{main} in which we drop some of the conditions.
Strictly speaking, these `corollaries' are really corollaries of the {\it proof} of Theorem \ref{main}.
The first of those cases is the classical Dedekind reciprocity law, Theorem \ref{dedreclaw} \cite{dedekind}. Recall that 
 \[ \s (a,b) = \frac{1}{4} \ \c \left( \begin{array}{c|cc} b & a & 1 \\
                                                           0 & 0 & 0 \\
                                                           0 & 0 & 0 \end{array} \right) \ . \] 
{\it Proof of Theorem} \ref{dedreclaw}. To modify our proof of Theorem \ref{main}, we have to consider the function
  \[ f(z) = \cot ( \pi a z ) \cot ( \pi b z ) \cot ( \pi z ) \ . \]
The residues are computed as above and yield the classical Dedekind sums;
the only difference is an additional pole of order three at $z=0$ (we may choose our
integration path $\gamma$ such that 0 is inside $\gamma$). Its residue is easily computed as
  \[ \Res \left( f(z) , z=0 \right) = - \frac{1}{3} \left( \frac{ a }{ b } + \frac{ 1 }{ ab } + \frac{ b }{ a } \right) \ .  \]
From here we can proceed as before.
\hfill {} $\Box$

It should be mentioned that a proof of Dedekind's reciprocity law along these lines is given in \cite[p.~21]{grosswald}.
In fact, it was this proof that motivated the proof of Theorem \ref{main}. 

The second special case is an identity equivalent to the reciprocity law for the Dedekind Bernoulli sums
\cite{apostol,carlitzbernoulli,mikolas}. 
The respective statement for them can be obtained through Corollary \ref{apocot}.
\begin{corollary} Let $ a_{0}, a_{1}, a_{2} \in \N $ be pairwise relatively prime, and $ m_{0}, m_{1}, m_{2} \in \N_{0} $
not all zero, such that $ m_{0} + m_{1} + m_{2} $ is even. Then
  \begin{eqnarray*} &\mbox{}& (-1)^{ m_{0} } m_{0}! \sum_{ { l_{1} , l_{2} \geq 0 } \atop { l_{1} + l_{2} = m_{0} } } \frac{ a_{1}^{ l_{1} } a_{2}^{ l_{2} } }{ l_{1}! l_{2}! }
       \ \c \left( \begin{array}{c|cc} a_{0} & a_{1} & a_{2} \\
                                       m_{0} & m_{1} + l_{1} & m_{2} + l_{2} \\
                                       0 & 0 & 0 \end{array} \right) \\ 
                    &\mbox{}& \qquad + (-1)^{ m_{1} } m_{1}! \sum_{ { l_{0} , l_{2} \geq 0 } \atop { l_{0} + l_{2} = m_{1} } } \frac{ a_{0}^{ l_{0} } a_{2}^{ l_{2} } }{ l_{0}! l_{2}! }
       \ \c \left( \begin{array}{c|cc} a_{1} & a_{0} & a_{2} \\
                                       m_{1} & m_{0} + l_{0} & m_{2} + l_{2} \\
                                       0 & 0 & 0 \end{array} \right) \\
                    &\mbox{}& \qquad + (-1)^{ m_{2} } m_{2}! \sum_{ { l_{0} , l_{1} \geq 0 } \atop { l_{0} + l_{1} = m_{2} } } \frac{ a_{0}^{ l_{0} } a_{1}^{ l_{1} } }{ l_{0}! l_{1}! }
       \ \c \left( \begin{array}{c|cc} a_{2} & a_{0} & a_{1} \\
                                       m_{2} & m_{0} + l_{0} & m_{1} + l_{1} \\
                                       0 & 0 & 0 \end{array} \right) \\
                    &\mbox{}& = \phi ( a_{0}, a_{1}, a_{2}, m_{0}, m_{1}, m_{2} ) \ , \end{eqnarray*}
where
  \begin{eqnarray*} &\mbox{}& \phi ( a_{0}, a_{1}, a_{2}, m_{0}, m_{1}, m_{2} ) := (-4)^{ (m_{0} + m_{1} + m_{2}) / 2 } \cdot \\
                    &\mbox{}& \qquad \left( \frac{ (-1)^{ m_{0} } }{ a_{0}^{ m_{0} + 1 } } \sum_{ { 2 k_{1} \geq m_{1} + 1 , \ 2 k_{2} \geq m_{2} + 1 } \atop { 2 ( k_{1} + k_{2} - 1 ) = m_{0} + m_{1} + m_{2} } } { m_{0} \choose { 2 k_{1} - 1 - m_{1} } } \frac{ B_{ 2 k_{1} } B_{ 2 k_{2} } }{ k_{1} k_{2} } a_{1}^{ 2 k_{1} - 1 - m_{1} } a_{2}^{ 2 k_{2} - 1 - m_{2} } \right. \\
                    &\mbox{}& \qquad + \frac{ (-1)^{ m_{1} } }{ a_{1}^{ m_{1} + 1 } } \sum_{ { 2 k_{0} \geq m_{0} + 1 , \ 2 k_{2} \geq m_{2} + 1 } \atop { 2 ( k_{0} + k_{2} - 1 ) = m_{0} + m_{1} + m_{2} } } { m_{1} \choose { 2 k_{2} - 1 - m_{2} } } \frac{ B_{ 2 k_{0} } B_{ 2 k_{2} } }{ k_{0} k_{2} } a_{0}^{ 2 k_{0} - 1 - m_{0} } a_{2}^{ 2 k_{2} - 1 - m_{2} } \\
                    &\mbox{}& \qquad + \left. \frac{ (-1)^{ m_{2} } }{ a_{2}^{ m_{2} + 1 } } \sum_{ { 2 k_{0} \geq m_{0} + 1 , \ 2 k_{1} \geq m_{1} + 1 } \atop { 2 ( k_{0} + k_{1} - 1 ) = m_{0} + m_{1} + m_{2} } } { m_{2} \choose { 2 k_{0} - 1 - m_{0} } } \frac{ B_{ 2 k_{0} } B_{ 2 k_{1} } }{ k_{0} k_{1} } a_{0}^{ 2 k_{0} - 1 - m_{0} } a_{1}^{ 2 k_{1} - 1 - m_{1} } \right) \\
                    &\mbox{}& \qquad + \frac{ 2^{ m_{0} + m_{1} + m_{2} + 2 } B_{ m_{0} + m_{1} + m_{2} + 2 } }{ m_{0} + m_{1} + m_{2} + 2 }
                                              \left( \frac{ (-1)^{ ( m_{0} + m_{1} - m_{2} ) / 2 } m_{0} ! m_{1} ! a_{2}^{ m_{0} + m_{1} + 1 } }{ ( m_{0} + m_{1} + 1 ) ! a_{0}^{ m_{0} + 1 } a_{1}^{ m_{1} + 1 } } \right. \\
                    &\mbox{}& \qquad \quad \left. + \frac{ (-1)^{ ( m_{0} + m_{2} - m_{1} ) / 2 } m_{0} ! m_{2} ! a_{1}^{ m_{0} + m_{2} + 1 } }{ ( m_{0} + m_{2} + 1 ) ! a_{0}^{ m_{0} + 1 } a_{2}^{ m_{2} + 1 } }
                                                  + \frac{ (-1)^{ ( m_{1} + m_{2} - m_{0} ) / 2 } m_{1} ! m_{2} ! a_{0}^{ m_{1} + m_{2} + 1 } }{ ( m_{1} + m_{2} + 1 ) ! a_{1}^{ m_{1} + 1 } a_{2}^{ m_{2} + 1 } } \right) \ . \end{eqnarray*}
\end{corollary}
Again the parity restriction on $ m_{0} + m_{1} + m_{2} $ is no constraint, since otherwise all the sums vanish.

{\it Proof.}
As in the last proof, we use the function 
  \[ f(z) = \cot^{ ( m_{0} ) } ( \pi a_{0} z ) \cot^{ ( m_{1} ) } ( \pi a_{1} z ) \cot^{ ( m_{2} ) } ( a_{2} \pi z ) \ . \]
By the pairwise-prime condition, the poles of $f$ are all simple with exception of the pole at $ z=0$ (again we may choose our
integration path such that 0 is contained inside). The residues at $ z = k/a_{0} $, for example, yield, as in the proof of Theorem \ref{main},
  \begin{eqnarray*} &\mbox{}& \sum_{ k=1 }^{ a_{0} - 1 } \Res \left( f(z) , z = \frac{ k }{ a_{0} } \right) 
                      = \frac{1}{\pi} (-1)^{ m_{0} } m_{0}! \sum_{ { l_{1} , l_{2} \geq 0 } \atop { l_{1} + l_{2} = m_{0} } } \frac{ a_{1}^{ l_{1} } a_{2}^{ l_{2} } }{ l_{1}! l_{2}! }
       \ \c \left( \begin{array}{c|cc} a_{0} & a_{1} & a_{2} \\
                                       m_{0} & m_{1} + l_{1} & m_{2} + l_{2} \\
                                       0 & 0 & 0 \end{array} \right) \ . \end{eqnarray*}
Similar expressions are obtained for the other non-zero poles of $f$. To get the residue at $ z=0 $, we use
the following expansion of the cotangent, which follows directly from the definition (\ref{berdef}) of the
Bernoulli polynomials \cite{bronstein}: 
  \begin{equation}\label{cotbern} \cot z = \frac{1}{z} + \sum_{ k \geq 1 } \frac{ (-1)^{k} 2^{ 2k } B_{2k} }{ (2k)! } z^{ 2k-1 } \ . \end{equation}
Hence
  \[ \cot^{(m)} z = \frac{ (-1)^{m} m! }{ z^{m+1} } + \sum_{ 2k \geq m+1 } \frac{ (-1)^{k} 2^{ 2k-1 } B_{2k} }{ k \ (2k-1-m)! } z^{ 2k-1-m } \ , \]
from which we obtain
  \[ \Res \left( f(z) , z = 0 \right) = - \frac{1}{\pi} \ \phi ( a_{0}, a_{1}, a_{2}, m_{0}, m_{1}, m_{2} ) \ . \]
Since not all $ m_{0}, m_{1}, m_{2} $ are zero, the sum of the residues of $f$ vanishes as in the proof of Theorem \ref{main},
and the statement follows.
\hfill {} $\Box$

Next we prove the reciprocity law for the cotangent sums \cite{dietercot}. Recall that
 \[ \c \left( a, b, c; x, y, z \right) = \ \c \left( \begin{array}{c|cc} c & a & b \\
                                                                         0 & 0 & 0 \\
                                                                         z & x & y \end{array} \right) \ . \]
It is easy to see that the reciprocity law of Meyer and Sczech \cite{sczech} is a special case of
\begin{corollary}[Dieter] Let $ a, b, c \in \N $ be pairwise relatively prime, and define $ A, B, C $ by
  \[ Abc + Bca + Cab = 1 \ . \]
Let $ x, y, z \in \R $, not all integers, and set
  \[ x' = cy - bz , y' = az - cx , z' = bx - ay \ . \]
Finally, let
  \[ \delta (x) = \left\{ \begin{array}{cl} 1 & \mbox{ if } x \in \Z , \\ 
                                             0 & \mbox{ else. } \end{array} \right. \]
Then
  \begin{eqnarray*} &\mbox{}& \c \left( a, b, c; x, y, z \right) + \c \left( b, c, a; y, z, x \right) + \c \left( c, a, b; z, x, y \right) = - 1 - \frac{ c }{ ab } \delta (z') \cot^{(1)} \pi \left( Acx' - (Bc+Cb) y' \right) \\
                    &\mbox{}& \qquad \quad - \frac{ a }{ bc } \delta (x') \cot^{(1)} \pi \left( Bay' - (Ca+Ac) z' \right) - \frac{ b }{ ac } \delta (y') \cot^{(1)} \pi \left( Cbz' - (Ab+Ba) x' \right) \ . \end{eqnarray*}
\end{corollary}
{\it Proof.} As in the proof of Theorem \ref{main}, we use the function
  \[ f(w) = \cot \pi ( aw - x ) \cot \pi ( bw - y ) \cot \pi ( cw - z ) \ . \]
If $f$ has only simple poles, the statement is a direct special case of Theorem \ref{main}. Otherwise, suppose we have
a double pole $w_{p}$, that is, there exist integers $m$ and $n$ such that (for example),
  \[ w_{p} = \frac{ m+x }{ a } = \frac{ n+y }{ b } \ . \]
To compute the residue of $f$ at this pole, we use once again the Laurent series (\ref{laurent}) and
(\ref{taylor}) of the cotangent to obtain
  \[ \Res \left( f(w) , w=w_{p} \right) = \frac{ c }{ \pi a b } \ \cot^{(1)} \pi \left( c w_{p} - z \right) \ . \]
We obtain similar residues for the simple poles as before, and the residue theorem yields an identity.
That this identity is equivalent to Dieter's can be easily seen by following the remarks in \cite{dietercot}
just before Theorem 2.3.
\hfill {} $\Box$

The last special case is one in general `dimension', Zagier's higher-dimensional Dedekind sums
\cite{zagier}. Recall that
 \[ \s \left( a_{0} ; a_{1}, \dots , a_{d} \right) = (-1)^{d/2} \ \c \left( \begin{array}{c|ccc} a_{0} & a_{1} & \cdots & a_{d} \\
                                                                                                 0 & 0 & \cdots & 0 \\
                                                                                                 0 & 0 & \cdots & 0 \end{array} \right) \ .  \]
\begin{corollary}[Zagier]\label{zagierthm} If $ a_{0} , \dots , a_{d} \in \N $ are pairwise relatively prime then
  \[ \sum_{ n=0 }^{ d } \s \left( a_{n} ; a_{0} , \dots , \widehat{ a_{n} } , \dots , a_{d} \right) = 1 - h ( a_{0} , \dots , a_{d} ) \ , \]
where
  \[ h ( a_{0} , \dots , a_{d} ) := \frac{ 2^{d} }{ a_{0} \cdots a_{d} } \sum_{ { k_{0} , \dots , k_{d} \geq 0 } \atop { k_{0} + \dots + k_{d} = d/2 } } \frac{ B_{ 2 k_{0} } \cdots B_{ 2 k_{d} }  }{ ( 2 k_{0} )! \cdots ( 2 k_{d} )!  } a_{0}^{ 2 k_{0} } \cdots a_{d}^{ 2 k_{d} } \ . \]
\end{corollary}
This description of $h ( a_{0} , \dots , a_{d} )$ in terms of Bernoulli numbers was first used by Berndt
\cite{berndt}; however, it is easily seen to be equivalent to the version given by Zagier. 
It is interesting to note that $h ( a_{0} , \dots , a_{d} )$ can be expressed in terms of Hirzebruch
L-functions \cite{zagier}. 

{\it Proof.} This time consider the function
  \[ f(z) = \cot \pi a_{0} z \ \cdots \ \cot \pi a_{d} z \ . \]
As before, $f$ has simple poles aside from the pole at $ z=0 $. The residues are calculated as usual, for example,
  \begin{eqnarray*} &\mbox{}& \sum_{ k=1 }^{ a_{0} - 1 } \Res \left( f(z) , z = \frac{ k }{ a_{0} } \right) = \frac{1}{ \pi a_{0} } \sum_{ k=1 }^{ a_{0} - 1 } \cot \frac{ \pi k a_{1} }{ a_{0} } \cdots \cot \frac{ \pi k a_{d} }{ a_{0} } = \frac{ (-1)^{ d/2 } }{ \pi } \ \s ( a_{0}; a_{1}, \dots , a_{d} ) \ .
\end{eqnarray*}
The residue at 0 can be computed through rewriting (\ref{cotbern}) as
  \[ z \cot z = \sum_{ k \geq 0 } \frac{ (-1)^{k} 2^{ 2k } B_{2k} }{ (2k)! } z^{ 2k } \ . \]
Hence
  \[ f(z) = \frac{ 1 }{ a_{0} \cdots a_{d} ( \pi z )^{ d+1 } } \prod_{ j=0 }^{ d } \sum_{ k_{j} \geq 0 } \frac{ (-1)^{ k_{j} } 2^{ 2k_{j} } B_{2k_{j}} }{ (2k_{j})! } ( \pi a_{j} z )^{ 2k_{j} } \ , \]
and we obtain the residue
  \begin{eqnarray*} &\mbox{}& \Res \left( f(z) , z = 0 \right) \frac{ (-1)^{ d/2 } 2^{d} }{ \pi a_{0} \cdots a_{d} } \sum_{ { k_{0} , \dots , k_{d} \geq 0 } \atop { k_{0} + \dots + k_{d} = d/2 } } \frac{ B_{ 2 k_{0} } \cdots B_{ 2 k_{d} }  }{ ( 2 k_{0} )! \cdots ( 2 k_{d} )!  } a_{0}^{ 2 k_{0} } \cdots a_{d}^{ 2 k_{d} } = \frac{ (-1)^{ d/2 } }{ \pi } \ h ( a_{0} , \dots , a_{d} ) \ .
\end{eqnarray*} 
It remains to apply the residue theorem and (\ref{residues}).
\hfill {} $\Box$


\section{Petersson-Knopp identities}\label{petkn}
Knopp applied in \cite{knopp} Hecke operators to $ \log \eta $ to arrive at Theorem \ref{petknopp}. 
This identity was stated by Petersson in the 1970's with additional congruence restrictions on $a$ and $b$.
For $n$ prime, the Petersson-Knopp identity was already known to Dedekind \cite{dedekind}. 
Theorem \ref{petknopp} was generalized by Parson and Rosen to Dedekind-Bernoulli sums \cite{rosen}, by
Apostol and Vu to their `sums of Dedekind type' \cite{vu}, 
and, most broadly, by Zheng to what we will call sums of Dedekind type with weight $ ( m_{1} , m_{2} ) $
\cite{zheng}. We will state and further generalize Zheng's Petersson-Knopp identity after the following
\begin{definition} Let $ a, a_{1}, \dots , a_{d} \in \N $. The sum 
  \[ S \left( a; a_{1}, \dots , a_{d} \right) := \sum_{ k \mbox{ \rm \scriptsize mod } a } f_{1} \left( \frac{k a_{1}  }{a } \right) \cdots f_{d} \left( \frac{k a_{d}  }{a } \right)  \]
is said to be {\bf of Dedekind type with weight $ \left( m_{1} , \dots , m_{d} \right) $} if for all
$ j = 1, \dots, d $, $ f_{j} (x+1) = f_{j} (x) $ and for all $ a \in \N$,
  \begin{equation}\label{weight} \sum_{ k \mbox{ \rm \scriptsize mod } a } f_{j} \left( x + \frac{ k }{ a } \right) = a^{ m_{j} } f_{j} (ax) \ . \end{equation}
\end{definition}
Note that the Bernoulli functions $ \B_{m} (x) $ satisfy (\ref{weight}) (with `weight' $ -m+1 $), as do the
functions $ \cot^{ (m) } ( \pi x ) $ (with `weight' $m+1$).
Zheng's theorem is the `two-dimensional' ($d=2$) case of the following 
\begin{theorem}\label{petknoppcot} Let $ n, a, a_{1}, \dots , a_{d} \in \N $. If
  \[ S \left( a; a_{1}, \dots , a_{d} \right) := \sum_{ k \mbox{ \rm \scriptsize mod } a } f_{1} \left( \frac{k a_{1}  }{a } \right) \cdots f_{d} \left( \frac{k a_{d}  }{a } \right)  \]
is of Dedekind type with weight $ \left( m_{1} , \dots , m_{d} \right) $ then
  \[ \sum_{ b | n } b^{ - m_{1} - \dots - m_{d} } \! \sum_{ r_{1} , \dots , r_{d} \mbox{ \rm \scriptsize mod } b } \! S \left( ab; \frac{ n }{ b } a_{1} + r_{1} a \ , \ \dots \ , \ \frac{ n }{ b } a_{d} + r_{d} a  \right) = \ n \ \sigma_{ d - 1 - m_{1} - \dots - m_{d}  } (n) \ S \left( a; a_{1}, \dots , a_{d} \right) \ . \] 
\end{theorem}
Our proof is a relatively straightforward extension of Zheng's proof for $d=2$ \cite{zheng}. 
We need the following two identities: 
\begin{lemma}\label{petknopplemma} Let $ a, a_{1}, \dots , a_{d} \in \N $. If
  \[ S \left( a; a_{1}, \dots , a_{d} \right) := \sum_{ k \mbox{ \rm \scriptsize mod } a } f_{1} \left( \frac{k a_{1}  }{a } \right) \cdots f_{d} \left( \frac{k a_{d}  }{a } \right)  \]
is of Dedekind type with weight $ \left( m_{1} , \dots , m_{d} \right) $ then for all $ j = 1, \dots, d $, 
  \begin{equation} \sum_{ k \mbox{ \rm \scriptsize mod } a } f_{j} \left( x + \frac{ k b }{ a } \right) = (a,b)^{ 1 - m_{j} } a^{ m_{j} } f_{j} \left( \frac{ ax }{ (a,b) } \right) \label{one} \end{equation}
and
  \begin{equation} S \left( a b; a_{1} b, \dots , a_{d} b \right) = b \ S \left( a; a_{1}, \dots , a_{d} \right) \ . \label{two} \end{equation}
\end{lemma}
{\it Proof.} (\ref{one}): If $ (a,b) = 1 $, the statement is vacuous. If $a$ and $b$ are not relatively
prime, let $ a' = a / (a,b) $ and $ b' = b / (a,b) $. Then
  \begin{align*} \sum_{ k \mbox{ \rm \scriptsize mod } a } f_{j} \left( x + \frac{ k b }{ a } \right) &= \sum_{ k \mbox{ \rm \scriptsize mod } a } f_{j} \left( x + \frac{ k b' }{ a' } \right) = (a,b) \sum_{ k \mbox{ \rm \scriptsize mod } a' } f_{j} \left( x + \frac{ k b' }{ a' } \right) \\
                                                                                                      &= (a,b) \sum_{ k \mbox{ \rm \scriptsize mod } a' } f_{j} \left( x + \frac{ k }{ a' } \right) = (a,b) {a'}^{m_{j}} f_{j} ( a' x ) = (a,b)^{ 1 - m_{j} } a^{m} f_{j} \left( \frac{ ax }{ (a,b) } \right) \ . \end{align*}
(\ref{two}):
  \begin{align*} S \left( a b; a_{1} b, \dots , a_{d} b \right) &= \sum_{ k \mbox{ \rm \scriptsize mod } ab } f_{1} \left( \frac{ k a_{1}  }{a } \right) \cdots f_{d} \left( \frac{k a_{d}  }{a } \right) \\
                                                                &= b \sum_{ k \mbox{ \rm \scriptsize mod } a } f_{1} \left( \frac{k a_{1}  }{a } \right) \cdots f_{d} \left( \frac{k a_{d}  }{a } \right) = b \ S \left( a; a_{1}, \dots , a_{d} \right) \ . \end{align*} 
\hfill {} $\Box$

{\it Proof of Theorem} \ref{petknoppcot}. We will make use of two properties of the M\"obius $ \mu $-function
  \[ \mu (n) = \left\{ \begin{array}{cl} 1 & \mbox{ if } n = 1 , \\ 
                                  (-1)^{m} & \mbox{ if } n = p_{1} \cdots p_{m} \mbox{ is square-free, } \\ 
                                         0 & \mbox{ otherwise, }  \end{array} \right. \]
namely,
  \begin{equation}\label{mu1} \sum_{ d|n } \mu (d) = \left\{ \begin{array}{cl} 1 & \mbox{ if } n = 1 , \\ 
                                                                               0 & \mbox{ otherwise, } \end{array} \right. \end{equation}
and
  \begin{equation}\label{mu} \sum_{ { k=1 } \atop { (k,b) = 1 } }^{ab} f(k) = \sum_{ t|b } \mu (t) \sum_{ k=1 }^{ ab/t } f(tk) \ .  \end{equation}
These suffice to prove our statement:
  \begin{eqnarray*} &\mbox{}& \sum_{ b | n } b^{ - m_{1} - \dots - m_{d} } \sum_{ r_{1} , \dots , r_{d} \mbox{ \rm \scriptsize mod } b } S \left( ab; \frac{ n a_{1} }{ b } + r_{1} a , \dots , \frac{ n a_{d} }{ b } + r_{d} a  \right) \\
                    &\mbox{}& \qquad = \sum_{ b | n } b^{ - m_{1} - \dots - m_{d} } \sum_{ { r_{1} , \dots , r_{d} \mbox{ \rm \tiny mod } b } \atop { \ k \mbox{ \rm \tiny mod } ab } } f_{1} \left( \frac{ k n a_{1} }{ a b^{2} } + \frac{ k r_{1} }{ b } \right) \cdots f_{d} \left( \frac{ k n a_{d} }{ a b^{2} } + \frac{ k r_{d} }{ b } \right) \\
                    &\mbox{}& \quad \ \ \stackrel{ \mbox{\rm \scriptsize (\ref{one})} }{ = } \sum_{ b | n } b^{ - m_{1} - \dots - m_{d} } \sum_{ k \mbox{ \rm \scriptsize mod } ab } (k,b)^{ 1 - m_{1} } b^{ m_{1} } f_{1} \left( \frac{ k n a_{1} }{ a b (a,b) } \right) \cdots (k,b)^{ 1 - m_{d} } b^{ m_{d} } f_{d} \left( \frac{ k n a_{d} }{ a b (a,b) } \right) \\
                    &\mbox{}& \qquad = \sum_{ b | n } \sum_{ k \mbox{ \rm \scriptsize mod } ab } (k,b)^{ d - m_{1} - \dots - m_{d} } f_{1} \left( \frac{ k n a_{1} }{ a b (a,b) } \right) \cdots f_{d} \left( \frac{ k n a_{d} }{ a b (a,b) } \right) \\
                    &\mbox{}& \qquad = \sum_{ b | n } \sum_{ c | b } c^{ d - m_{1} - \dots - m_{d} } \sum_{ { k \mbox{ \rm \tiny mod } ab } \atop { (k,b) = c } } f_{1} \left( \frac{ k n a_{1} }{ a b c } \right) \cdots f_{d} \left( \frac{ k n a_{d} }{ a b c } \right) \\
                    &\mbox{}& \qquad = \sum_{ b | n } \sum_{ c | b } c^{ d - m_{1} - \dots - m_{d} } \sum_{ { k \mbox{ \rm \tiny mod } ab/c } \atop { (k,b/c) = 1 } } f_{1} \left( \frac{ k n a_{1} }{ a b } \right) \cdots f_{d} \left( \frac{ k n a_{d} }{ a b } \right) \\
                    &\mbox{}& \quad \ \ \stackrel{ \mbox{\rm \scriptsize (\ref{mu})} }{ = } \sum_{ b | n } \sum_{ c | b } c^{ d - m_{1} - \dots - m_{d} } \sum_{ t | b/c } \mu (t) \sum_{ k \mbox{ \rm \scriptsize mod } ab/ct } f_{1} \left( \frac{ t k n a_{1} }{ a b } \right) \cdots f_{d} \left( \frac{ t k n a_{d} }{ a b } \right) \\
                    &\mbox{}& \qquad = \sum_{ b | n } \sum_{ c | b } c^{ d - m_{1} - \dots - m_{d} } \sum_{ t | b/c } \mu (t) \ S \left( \frac{ ab }{ ct } ; \frac{ n a_{1} }{ c } , \dots , \frac{ n a_{d} }{ c } \right) \\
                    &\mbox{}& \qquad = \sum_{ cte | n } c^{ d - m_{1} - \dots - m_{d} } \mu (t) \ S \left( ae ; \frac{ n a_{1} }{ c } , \dots , \frac{ n a_{d} }{ c } \right) \\
                    &\mbox{}& \qquad = \sum_{ ce | n } c^{ d - m_{1} - \dots - m_{d} } \ S \left( ae ; \frac{ n a_{1} }{ c } , \dots , \frac{ n a_{d} }{ c } \right) \sum_{ t | n / ce } \mu (t) \\
                    &\mbox{}& \quad \ \ \stackrel{ \mbox{\rm \scriptsize (\ref{mu1})} }{ = } \sum_{ ce = n } c^{ d - m_{1} - \dots - m_{d} } \ S \left( ae ; a_{1} e , \dots , a_{d} e \right) \\
                    &\mbox{}& \quad \ \ \stackrel{ \mbox{\rm \scriptsize (\ref{two})} }{ = } n \ \sigma_{ d - 1 - m_{1} - \dots - m_{d}  } (n) \ S \left( a; a_{1}, \dots , a_{d} \right) \ . \end{eqnarray*}
\hfill {} $\Box$

The fact that the Dedekind cotangent sums
  \[ a_{0}^{ m_{0} + 1 } \ \c \left( \begin{array}{c|ccc} a_{0} & a_{1} & \cdots & a_{d} \\
                                                          m_{0} & m_{1} & \cdots & m_{d} \\
                                                              0 & 0 & \cdots & 0 \end{array} \right)
        = \sum_{ k \mbox{ \rm \scriptsize mod } a_{0} } \prod_{ j=1 }^{ d } \cot^{ ( m_{j} ) } \frac{ \pi k a_{j} }{ a_{ 0 } } \]
have weight $ ( m_{1} + 1 , \dots , m_{d} + 1 ) $ immediately yields the Petersson-Knopp-like Theorem \ref{petknoppcotangentsum}. 

A particularly simple form of Theorem \ref{petknoppcotangentsum} is achieved for Zagier's higher-dimensional 
Dedekind sums (the case $ m_{0} = \dots = m_{d} = 0 $): with
 \[ n \ \sigma_{ -1 } (n) = \sum_{ d|n } \frac{n}{d} = \sum_{ d|n } d = \sigma (n) \ , \]
we obtain 
\begin{corollary} For $ n, a_{0}, \dots , a_{d} \in \N $,
  \[ \sum_{ b | n } b^{ 1-d } \sum_{ r_{1} , \dots  ,r_{d} \mbox{ \rm \scriptsize mod } b } \s ( a_{0} b ; \frac{ n }{ b } a_{1} + r_{1} a_{0} , \dots , \frac{ n }{ b } a_{d} + r_{d} a_{0} ) = \sigma (n) \ \s ( a_{0} ; a_{1} , \dots , a_{d} ) \ . \] 
\hfill {} $\Box$
\end{corollary}


\section{Proof of the polynomial-time computability}\label{comp}
In this last section, we prove that the Dedekind cotangent sums are computable in polynomial time (Theorem \ref{compcotangent}). 
Our proof is similar to the last section of \cite{br}, where the computability of Zagier's
higher-dimensional Dedekind sums was shown. Again we will `merge' two theorems in combinatorial 
geometry, due to Barvinok \cite{barvinok} and Diaz-Robins \cite{diaz}. The latter
allows us to express generating functions for the lattice point count in cones in terms of
cotangents; Barvinok's theorem tells us that the respective rational functions are computable.
The usual generating function for lattice point enumeration for a $d$-dimensional cone $\K$, 
  \[ F_{\K} ( \q ) = F_{\K} ( q_{1} , \dots , q_{d} ) := \sum_{ \m = ( m_{1} , \dots , m_{d} ) \in \K \cap \Z^{d} } q_{1}^{ m_{1} } \cdots q_{d}^{ m_{d} } \ , \] 
assumes in \cite{diaz} an exponential variable:
\begin{theorem}[Diaz-Robins]\label{sinaithm} Suppose the cone $ \K \subset \R^{d} $ is generated by the positive
real span of the integer vectors $ \v_{1} , \dots , \v_{d} \in \Z^{d} $ such that the 
$d \times d $-matrix $ M = ( a_{ij} ) $, whose column vectors are $ \v_{1} , \dots , \v_{d} $, is lower-triangular. Let
$ p_{k} := a_{11} \cdots a_{kk} \ ( k = 1, \dots , d ) $ and
$ G := \left( \Z / p_{1} \Z \right) \times \dots \times \left( \Z / p_{d} \Z \right) $. Then
  \[ \sum_{ \m \in \K \cap \Z^{d} } e^{ - 2 \pi \left< \m , \cs \right> }  = \frac{ 1 }{ 2^{d} |G| } \sum_{ \r \in G } \prod_{ k=1 }^{ d } \left( 1 + \coth \frac{ \pi }{ p_{k} } \left< \cs + i \r , \v_{k} \right> \right) \ . \]
Here $ \left< \ , \ \right> $ denotes the usual scalar product in $ \R^{d} $.
\end{theorem} 
We note that the assumption on $M$ being lower-triangular is not crucial for practical purposes: by Hermite
normal form, any cone generated by integer vectors is unimodular equivalent to a cone described by a lower-triangular matrix.
\begin{theorem}[Barvinok]\label{barvthm} For fixed dimension $d$, the rational function $ F_{ \K } ( \q ) $
is polynomial-time computable in the input size of $\K$.
\end{theorem}
These two powerful theorems combined allow us to compute the Dedekind cotangent sum in polynomial time: 

{\it Proof of Theorem} \ref{compcotangent}. 
It suffices to prove the polynomial-time computability of
  \[ \sum_{ k \mbox{ \rm \scriptsize mod } a } \ \prod_{ j=1 }^{ d } \cot \pi \left( \frac{ k a_{j} }{ a } + z_{j} \right) \ . \]
We will do this inductively by constructing a cone whose generating function has the above Dedekind cotangent
sum + `lower-dimensional' sums, that is, Dedekind cotangent sums with a lower number of factors.
More precisely, we will construct a cone whose generating function will be
  \[ \sum_{ k \mbox{ \rm \scriptsize mod } a } \left( 1 + \cot \pi \left( \frac{ k a_{1} }{ a } + z_{1} \right) \right) \cdots \left( 1 + \cot \pi \left( \frac{ k a_{d} }{ a } + z_{d} \right) \right) \ . \]
By induction and Barvinok's Theorem \ref{barvthm}, this will prove the computability of the Dedekind cotangent sum.
Let $\K$ be the positive real span of 
  \[ \left( \begin{array}{c} a \\ 0 \\ \vdots \\ \\ 0 \\ b_{1} \end{array} \right) ,
     \left( \begin{array}{c} 0 \\ a \\ 0 \\ \vdots \\ 0 \\ b_{2} \end{array} \right) , \dots ,
     \left( \begin{array}{c} 0 \\ \vdots \\ \\ 0 \\ a \\ b_{d-1} \end{array} \right) ,
     \left( \begin{array}{c} 0 \\ \\ \vdots \\ \\ 0 \\ a \end{array} \right) \ . \]
The integers $ b_{1} , \dots , b_{d-1} $ are to be chosen later. We will repeatedly use the fact that
  \begin{equation}\label{work} \sum_{ k \mbox{ \rm \scriptsize mod } a } \coth \pi \left( \frac{ik}{a} + z \right) = a \coth \pi a z \ . \end{equation}
To apply Theorem \ref{sinaithm} to $\K$, note that in our case $ p_{k} = a^{k} $. Hence 
  \begin{eqnarray*} &\mbox{}& \sum_{ \m \in \K \cap \Z^{d} } e^{ - 2 \pi \left< \m , \cs \right> }  = \frac{ 1 }{ 2^{d} a^{ d(d+1)/2 } } \sum_{ {r_{j} \mbox{ \rm \tiny mod } a^{j} } \atop { ( j = 1, \dots , d ) } } \prod_{ k=1 }^{ d } \left( 1 + \coth \frac{ \pi }{ a^{k} } \left< \cs + i \r , \v_{k} \right> \right) \\
                    &\mbox{}& \qquad = \frac{ 1 }{ 2^{d} a^{ d(d+1)/2 } } \sum_{ {r_{j} \mbox{ \rm \tiny mod } a^{j} } \atop { ( j = 1, \dots , d ) } } \left( 1 + \coth \frac{ \pi }{ a } \left( ( s_{1} + i r_{1} ) a + ( s_{d} + i r_{d} ) b_{1} \right) \right) \\
                    &\mbox{}& \qquad \qquad \left( 1 + \coth \frac{ \pi }{ a^{2} } \left( ( s_{2} + i r_{2} ) a + ( s_{d} + i r_{d} ) b_{2} \right) \right) \cdots \left( 1 + \coth \frac{ \pi }{ a^{d} } ( s_{d} + i r_{d} ) a \right)  \\
                    &\mbox{}& \qquad = \frac{ 1 }{ 2^{d} a^{ d(d+1)/2 - 1 } } \sum_{ {r_{j} \mbox{ \rm \tiny mod } a^{j} } \atop { ( j = 3, \dots , d ) } } \left( 1 + \coth \pi \left( s_{1} + ( s_{d} + i r_{d} ) \frac{ b_{1} }{ a } \right) \right) \\
                    &\mbox{}& \qquad \qquad \left( 1 + \coth \frac{ \pi }{ a^{3} } \left( ( s_{3} + i r_{3} ) a + ( s_{d} + i r_{d} ) b_{3} \right) \right) \cdots \left( 1 + \coth \frac{ \pi }{ a^{d-1} } ( s_{d} + i r_{d} ) \right) \\
                    &\mbox{}& \qquad \qquad \sum_{ r_{2} \mbox{ \rm \scriptsize mod } a^{2} } \left( 1 + \coth \pi \left( \frac{ i r_{2} }{ a } + \frac{ s_{2} }{ a } + ( s_{d} + i r_{d} ) \frac{ b_{2} }{ a^{2} } \right) \right) \\
                    &\mbox{}& \quad \ \ \stackrel{ \mbox{\rm \scriptsize (\ref{work})} }{ = } \frac{ 1 }{ 2^{d} a^{ d(d+1)/2 - 3 } } \sum_{ {r_{j} \mbox{ \rm \tiny mod } a^{j} } \atop { ( j = 4, \dots , d ) } } \left( 1 + \coth \pi \left( s_{1} + ( s_{d} + i r_{d} ) \frac{ b_{1} }{ a } \right) \right) \\
                    &\mbox{}& \qquad \qquad \left( 1 + \coth \pi \left( s_{2} + ( s_{d} + i r_{d} ) \frac{ b_{2} }{ a } \right) \right) \left( 1 + \coth \frac{ \pi }{ a^{4} } \left( ( s_{4} + i r_{4} ) a + ( s_{d} + i r_{d} ) b_{4} \right) \right) \\
                    &\mbox{}& \qquad \qquad \cdots \left( 1 + \coth \frac{ \pi }{ a^{d-1} } ( s_{d} + i r_{d} ) \right) \sum_{ r_{3} \mbox{ \rm \scriptsize mod } a^{3} } \left( 1 + \coth \pi \left( \frac{ i r_{3} }{ a^{2} } + \frac{ s_{3} }{ a^{2} } + ( s_{d} + i r_{d} ) \frac{ b_{3} }{ a^{3} } \right) \right) \\
                    &\mbox{}& \quad \ \ \stackrel{ \mbox{\rm \scriptsize (\ref{work})} }{ = } \dots = \frac{ 1 }{ 2^{d} a^{ d } } \sum_{ {r_{d} \mbox{ \rm \scriptsize mod } a^{d} } } \left( 1 + \coth \pi \left( s_{1} + ( s_{d} + i r_{d} ) \frac{ b_{1} }{ a } \right) \right) \cdots \\
                    &\mbox{}& \qquad \qquad \left( 1 + \coth \pi \left( s_{d-1} + ( s_{d} + i r_{d} ) \frac{ b_{d-1} }{ a } \right) \right) \left( 1 + \coth \frac{ \pi }{ a^{d-1} } ( s_{d} + i r_{d} ) \right)  \\
                    &\mbox{}& \qquad = \frac{ 1 }{ 2^{d} a^{ d } } \sum_{ n=1 }^{ a^{ d-1 } } \sum_{ k=1 }^{ a } \left( 1 + \coth \pi \left( s_{1} + \left( s_{d} + i (na+k) \right) \frac{ b_{1} }{ a } \right) \right) \cdots \\
                    &\mbox{}& \qquad \qquad \left( 1 + \coth \pi \left( s_{d-1} + \left( s_{d} + i (na+k) \right) \frac{ b_{d-1} }{ a } \right) \right) \left( 1 + \coth \frac{ \pi }{ a^{d-1} } \left( s_{d} + i (na+k) \right) \right)  \\
                    &\mbox{}& \qquad = \frac{ 1 }{ 2^{d} a^{ d-1 } } \sum_{ k=1 }^{ a } \left( 1 + \coth \pi \left( s_{1} + ( s_{d} + i k ) \frac{ b_{1} }{ a } \right) \right) \cdots \\
                    &\mbox{}& \qquad \qquad \left( 1 + \coth \pi \left( s_{d-1} + ( s_{d} + i k ) \frac{ b_{d-1} }{ a } \right) \right) \sum_{ n=1 }^{ a^{ d-2 } } \left( 1 + \coth \pi \left( \frac{ in }{ a^{d-2} } + \frac{ s_{d} }{ a^{d-1} } + \frac{ ik }{ a^{d-1} } \right) \right)  \\
                    &\mbox{}& \quad \ \ \stackrel{ \mbox{\rm \scriptsize (\ref{work})} }{ = } \frac{ 1 }{ 2^{d} a } \sum_{ k=1 }^{ a } \left( 1 + \coth \pi \left( s_{1} + ( s_{d} + i k ) \frac{ b_{1} }{ a } \right) \right) \cdots \\
                    &\mbox{}& \qquad \qquad \left( 1 + \coth \pi \left( s_{d-1} + ( s_{d} + i k ) \frac{ b_{d-1} }{ a } \right) \right) \left( 1 + \coth \pi \left( \frac{ s_{d} }{ a } + \frac{ ik }{ a } \right) \right)  \ . \end{eqnarray*} 
If we now choose $ b_{j} = a_{d}^{-1} a_{j} $ ($ j = 1 , \dots , d-1 $), where $ a_{d}^{-1} a_{d} \equiv 1 $ mod $a$, and
  \[ s_{j} = i z_{j} - i z_{d} b_{j} \ \ ( j = 1, \dots , d-1 ) , \qquad s_{d} = i a z_{d} \ ,  \]
this generating function becomes
  \[ \sum_{ \m \in \K \cap \Z^{d} } e^{ - 2 \pi \left< \m , \cs \right> } = \frac{ 1 }{ 2^{d} a } \sum_{ k \mbox{ \rm \scriptsize mod } a } \left( 1 + \cot \pi \left( \frac{ k a_{1} }{ a } + z_{1} \right) \right) \cdots \left( 1 + \cot \pi \left( \frac{ k a_{d} }{ a } + z_{d} \right) \right) \ . \] 
\hfill {} $\Box$


\section{Closing remarks} 
There remain many open questions. First, for all Dedekind sums there are closed formulas for special cases, for example, $ \s (1,a) $. 
Many of such formulas can be found in \cite{bruce}. The Dedekind cotangent sum and its variations attain such closed formulas for special 
variables; it is not clear how far these cases lead. Second, it seems desirable to have a parallel theory for the respective Dedekind sum 
in which the cotangents get replaced by Bernoulli functions. This will most certainly require different methods than the ones used in this 
paper. Finally, we can apply our very general Petersson-Knopp-like Theorem \ref{petknoppcot} to various other sums, for example, to such 
'Dedekind Bernoulli sums'. 

{\bf Acknowledgements}. I am grateful to Sinai Robins for good discussions and 
Thomas Zaslavsky for many helpful comments on a previous version of this paper. 


\bibliographystyle{alpha}

 \sc Department of Mathematical Sciences\\
 State University of New York\\
 Binghamton, NY 13902-6000\\
 {\tt matthias@math.binghamton.edu} \\ 
 {\tt http://www.math.binghamton.edu/matthias} 

\end{document}